%% file: SuperLRrule.tex
\documentclass[a4paper, twoside, 11pt, english]{article}

\usepackage[T1]{fontenc}
\usepackage[utf8]{inputenc}
\usepackage{etoolbox}
\usepackage{amsmath,amsthm,amssymb,stmaryrd}
\usepackage{mathtools}
\usepackage[ttscale=.875]{libertine}
\usepackage[libertine]{newtxmath}
\usepackage[numbers,sort&compress]{natbib}
\usepackage{comment}
\usepackage{ifthen}
\usepackage{algpseudocode}
\usepackage{algorithm2e}
\usepackage{multicol}
\usepackage{fancyhdr, titlesec, url, enumerate, microtype,setspace}
\usepackage[all,knot,poly]{xy}
\usepackage{tikz}
\usetikzlibrary{decorations.pathreplacing}

\usepackage{tikz-qtree,varwidth}
\usepackage{youngtab}
\usepackage{ytableau}
\usepackage[english]{babel}
\usepackage{hyperref}

\usetikzlibrary{calc}
\CompileMatrices
\input{macros}

\begin{document}
\thispagestyle{empty}

\begin{center}

\begin{doublespace}
\begin{huge}
{\scshape A super Littlewood--Richardson type rule}
\end{huge}

\bigskip
\hrule height 1.5pt 
\bigskip

\begin{Large}
{\scshape Nohra Hage }
\end{Large}

\end{doublespace}

\vspace{1cm} 

\begin{small}\begin{minipage}{14cm}
\noindent\textbf{Abstract -- }
We introduce a super version of the Littlewood--Richardson rule for super Schur functions over signed alphabets. 
We give in particular combinatorial interpretations of the super Littlewood--Richardson coefficients using the properties of super Young tableaux, which have found rich applications in representation theory, algebraic combinatorics, and mathematical physics.

\medskip


\medskip

\smallskip
\noindent
\textbf{M.S.C. 2010 -- Primary:} 05E10, 17B99.
\textbf{Secondary:} 05A99,  20M05.
\end{minipage}\end{small}
\end{center}

\tikzset{every tree node/.style={minimum width=1em,draw,circle},
         blank/.style={draw=none},
         edge from parent/.style=
         {draw,edge from parent path={(\tikzparentnode) -- (\tikzchildnode)}},
         level distance=0.8cm}
\bigskip

\section{Introduction}
The super algebraic structures have been used
as combinatorial tools in the study of the invariant theory and  the representation
theory of superalgebras, and algebras satisfying identities, and have played an important role in various areas of theoretical and mathematical physics,~\cite{BereleRemmel85,  BonettiSenatoVenezia98,GrosshansRotaStein87, KacWakimoto88}.

We are interested on the combinatorial study of super Young tableaux over signed alphabets of any order. These objects form a set of representatives of elements of the super plactic monoid, which is defined independently in~\cite{BenkartKangKashiwara20, LaScalaNardozzaSenato06,LodayPovov08} as a generalization of the classical plactic monoid over signed alphabets.
The plactic monoid has been originally introduced in~\cite{Schensted61, Knuth70, LascouxSchutsenberger81}, and has found several applications in representation theory, algebraic and probabilistic combinatorics, and rewriting theory,~\cite{LascouxSchutsenberger81, Hage15, HageMalbos17}. 
Super versions of the Robinson--Schensted--Knuth correspondence and the Littlewood--Richardson rule have been also introduced  over signed alphabets,~\cite{BereleRemmel85, BonettiSenatoVenezia98, KangKwon01, LaScalaNardozzaSenato06, Muth19, Hage2021Super, Hage2022RSK}. 
We introduce in this paper a  super Robinson--Schensted correspondence over signed alphabets  and  we present a super version of the Littlewood--Richardson rule in terms of the classical Littlewood--Richardson tableaux.

\section{A super Robinson--Schensted type correspondence}
\label{SS:SuperRobinsonSchenstedTypeCorrespondence}

We recall the notion of super tableaux  from~\cite{LaScalaNardozzaSenato06, Hage2021Super, Hage2022RSK}, and we introduce a super Robinson--Schensted correspondence over signed alphabets.

Denote $[n]:=\{1<\ldots<n\}$, for $n>0$. 
A \emph{signed alphabet} is a pair~$(\Si, ||.||)$ where~$\Si$ is a finite ordered set and~$||.||$ is any map from~$\Si$ to  the additive cyclic group~$\mathbb{Z}_2 = \{0, 1\}$ of order~$2$.  Denote~$\Si_0$ (resp.~$\Si_1$) the set of~$x$ in~$\Si$ such that~$||x||=0$                 (resp.~$||x||=1$).
A monoid~$\M$ is a \emph{supermonoid} if a map~$||.||:\M\to\mathbb{Z}_2$ is given such that~$||uv|| = ||u||+||v||$, for all~$u, v$ in~$\M$. Denote~$\Si^\ast$ the free monoid of words over~$\Si$.
A \emph{partition} of $n>0$ of \emph{height}~$k$ is a weakly decreasing
sequence~$\lambda=(\lambda_1,\ldots,\lambda_k)\in\mathbb{N}^k$ such that~$\sum \lambda_i = n$.
Denote~$\Pr_n$ the set of partitions of~$n$ and set~$P =\bigcup\Pr_n$. 
The \emph{Young diagram}
of a partition $\lambda= (\lambda_1,\ldots,\lambda_k)$, denoted by~$\Yr(\lambda)$, is the set of pairs~$(i,j)$ such that~$1\leq i\leq k$ and~$1\leq j\leq \lambda_i$, which is  represented by a diagram by drawing a box for each of its pairs.
A \emph{super tableau of shape}~$\lambda$ over~$\Si$, is a
pair~$t:=~(\lambda, \Tr)$, where~$\Tr:\Yr(\lambda)\to\Si$ is a map satisfying:
\begin{enumerate}[\bf C1)]
\item $\Tr(i,j)\leq~\Tr(i,j+1)$ with~$\Tr(i,j)=~\Tr(i,j+1)$ only if~$\Tr(i,j)\in\Si_0$, 
\item $\Tr(i,j)\leq \Tr(i+1,j)$ with~$\Tr(i,j)= \Tr(i+1,j)$ only if~$\Tr(i,j)\in\Si_1$. 
\end{enumerate}
The set~$\Yr(\lambda)$ (resp. map~$\Tr$) is called the \emph{frame} (resp.  \emph{filling}) of~$t$.
When~$\Tr$ is injective,~$t$ is called a \emph{standard tableau} over~$\Si$.
Denote by~$\YoungST(\Si)$ the set of all standard tableaux over~$\Si$,
and by~$\YoungT(\Si)$ (resp.~$\YoungT(\Si,\lambda)$) the set of all super tableaux (resp. of shape~$\lambda$) over~$\Si$.
An \emph{empty tableau of shape}~$\lambda$ is a super tableau of shape~$\lambda$ where all its boxes are empty. If there is no confusion, suppose that all empty tableaux are equal and denote by~$\emptyset$ the corresponding empty  tableau.
Let $\lambda, \mu\in\Pr$ of heights $k$ and $l$, respectively, such that~$l\leq k$. Denote by~$\mu\subseteq\lambda$ if~$\mu_i\leq \lambda_i$, for any~$i$,  and by~$\Yr(\lambda/\mu)$ the set of pairs~$(i,j)$ such that~$1\leq i\leq k$ and~$\mu_i< j\leq~\lambda_i$.
A \emph{super skew
tableau}  of \emph{shape} $\lambda/\mu$ over~$\Si$, is a pair $S := (\lambda/\mu, \U)$ where $\U: \Yr(\lambda/\mu)\to\Si$  is a map satisfying~${\bf C1)}$ and~${\bf C2)}$.
The set~$\Yr(\lambda/\mu)$ (resp.  map~$\U$) is called the \emph{frame} (resp. \emph{filling}) of~$S$.
We call~$\lambda/\mu$ a \emph{skew shape} and set $\lambda/0:= \lambda$.
Denote by~$\SkewT(\Si)$ (resp.~$\SkewT(\Si,\lambda/\mu)$) the set of all super skew tableaux (resp. of shape~$\lambda/\mu$) over~$\Si$. 
An \emph{inner corner} of~$S\in\SkewT(\Si,\lambda/\mu)$ is a box in~$\Yr(\mu)$ such that the
boxes below and to the right are not in~$\Yr(\mu)$, and 
an \emph{outer corner}  is a box such that neither box
below or to the right is in~$\Yr(\lambda/\mu)$. 
Denote by
\[
\r: \SkewT(\Si)\to\Si^\ast \quad \big(\text{resp}.~\c: \SkewT(\Si)\to\Si^\ast\big)
\]
the \emph{row} (resp. \emph{column}) \emph{reading} that reads a  tableau row (resp. column)-wise  from bottom to top and from left to right.

The \emph{column insertion}, denoted by $\insl{}$, inserts an element~$x$ in~$\Si$ into a super tableau~$t$ of~$\YoungT(\Si)$ as follows,~\cite{LaScalaNardozzaSenato06}.
If $t$ is empty, create a box and label it~$x$. 
Suppose $t$ is non-empty. If~$x\in\Si_0$ (resp.~$x\in\Si_1$) is larger than (resp. at least as large as) the bottom element of the leftmost column of~$t$, then put~$x$ in a box to the bottom of this column; Otherwise, let~$y$ be the smallest element of this column such that~$y\geq x$ (resp.~$y>x$). Then replace~$y$ by $x$ and recursively insert~$y$ into the super tableau formed by the columns of~$t$ to the right of the leftmost.
Note that this recursion may end with an insertion into an empty column to the right of the
existing columns of~$t$. Output the resulting super tableau as~$(x\insl{}t)$.
One defines in a similar way the \emph{row insertion}, denoted by $\insr{}$, that inserts an element in~$\Si$ into the rows of a super tableau starting from its topmost one,~\cite{LaScalaNardozzaSenato06}.
Define
\[
\T(w)  :=  (x_1\insl{}(\ldots\insl{}(x_k\insl{}\emptyset)\ldots)),
\]
for any~$w=x_1\ldots x_k$ in~$\Si^\ast$.

The \emph{super plactic monoid} over~$\Si$, denoted by~$\P(\Si)$, is the quotient of~$\Si^\ast$
by the  family of \emph{super Knuth relations},~\cite{BenkartKangKashiwara20, LaScalaNardozzaSenato06,LodayPovov08}:
\[
\begin{array}{rl}
xzy=zxy,&  \text{ with } x=y  \text{ only if }  ||y||=0\;  \text{ and }  y=z  \text{ only if }   ||y||=1,\\
yxz=yzx,& \text{ with } x=y  \text{ only if }  ||y||=1\;  \text{ and } y=z   \text{ only if } ||y||=0, 
\end{array}
\]
for any~$x\leq y\leq z\in\Si$.
Denote by~$\sim_{\P(\Si)}$ the super plactic congruence.
The set~$\YoungT(\Si)$ satisfies the \emph{cross-section property}  for~$\sim_{\P(\Si)}$:~$w\sim_{\P(\Si)} w'$  if and only if $\T(w)=\T(w')$, for all~$ w,w'$ in~$\Si^\ast$,~\cite{LaScalaNardozzaSenato06}.  
In particular, we have
\[
(y\insl{}(t\insr{}x)) = ((y\insl{}t)\insr{}x),\]
for all~$t$ in~$\YoungT(\Si)$ and~$x,y$ in~$\Si$. Moreover, we have~$\T(w)=((\ldots(\emptyset\insr{}x_1)\ldots)\insr{}x_k)$,  for any~$w=x_1\ldots x_k$ in~$\Si^\ast$. 

Define the following products on $\YoungT(\Si)$:
\[
\begin{array}{rl}
(t\insr{r}t') = (t\insr{}\r(t'))& \quad \big(\text{resp. }  (t'\insl{r}t) =  (\r(t')\insl{}t)\big),\\
(t\insr{c}t')= (t\insr{}\c(t'))&\quad  \big(\text{resp. }  (t'\insl{c}t) = (\c(t')\insl{}t)\big),
 \end{array}
 \]
 for all $t, t'$ in $\YoungT(\Si)$.
By the cross-section property, we have  
\[
(t\insr{r}t')= (t\insl{r}t')= (t\insr{c}t')= (t\insl{c}t'),
\]
\[
((t\insr{r}t')\insr{r}t'')=(t\insr{r}(t'\insl{r}t''))  = ((t\insl{r}t')\insl{r}t''),
\]
and 
\[
((t\insr{c}t')\insr{c}t'')=(t\insr{c}(t'\insr{c}t'')) = ((t\insl{c}t')\insl{c}t''),
\]
for all~$t, t', t''$ in $\YoungT(\Si)$.

Starting from a word $w=x_1\ldots x_k$ over $\Si$, we compute a pair of same shape tableaux $(\T(w),\Q(w))\in \YoungT(\Si)\times\YoungST([n])$, as follows.
Start with  empty super tableaux~$\T_{k+1}$ and~$\Q_{k+1}$. For each~$i=k,\ldots,1$,  compute~$x_i\insl{}\T_{i+1}$; let~$\T_i$ be the resulting super tableau. Add a
box filled with~($k-i+1$) to~$\Q_{i+1}$ in the same place as the box that belongs to~$\T_i$ but not
to~$\T_{i+1}$; let~$\Q_i$ be the resulting standard tableau. Output~$\T_1$ for $\T(w)$ and~$\Q_1$ as $\Q(w)$.

Let~$t$ be in~$\YoungT(\Si)$ and~$x$ in~$\Si$. 
Starting from~$(x\insl{}t)$, together with the outer corner that has been added to the frame of~$t$ after insertion, we  recover~$t$ and~$x$ as follows.
Suppose that $y\in \Si_1$ (resp.~$y\in\Si_0$) is the entry in the outer corner of~$t$. Find in the column to the right of this outer corner the largest entry which is strictly smaller (resp. smaller)
than~$y$. Then this entry is replaced by~$y$ and it is bumped right to the next column where the process is repeated until an entry is bumped out of the leftmost column. Output the resulting super tableau for~$t$ and the last element which is bumped out for~$x$.

Starting from a pair of same shape tableaux $(\T,\Q)\in\YoungT(\Si)\times\YoungST([n])$  consiting of $k$ boxes, we compute a word~$w\in\Si^\ast$ as follows.
Start with $\T_k = \T$ and $\Q_k = \Q$. For each~$i=k-1,\ldots, 1$,  take the box filled with $i+1$ in~$\Q_{i+1}$, and apply the reverse insertion algorithm to~$\T_{i+1}$ with the outer corner corresponding to that box. The
resulting tableau is~$\T_i$, and the element that is bumped out of the leftmost column of~$\T_{i+1}$ is denoted by~$x_{i+1}$.
Remove the box containing~$i+1$ from~$\Q_{i+1}$ and denote the resulting tableau by~$\Q_i$. Output the word~$x_1\ldots x_k$ as~$w$.

As a consequence, we deduce the following result:

\begin{theorem}
The map~$w\mapsto (\T(w),\Q(w))$ defines a bijection between words over~$\Si$ and
pairs of same-shape tableaux in~$\YoungT(\Si)\times\YoungST([n])$.
\end{theorem}

\section{A super Littlewood--Richardson type rule}
\label{S:ASuperLittlewoodRichardsonTypeRule}

Following the approach developed in~\cite{Thomas78} for the non-signed case,  we introduce a super version of the Littlewood--Richardson rule over signed alphabets.

Consider~$\lambda, \mu, \nu\in\Pr$ such that the partition~$\lambda\subseteq\nu$ is of height~$h$. Consider~$t\in\YoungT(\Si,\lambda)$ and~$t'\in\YoungT(\Si,\mu)$ such that~$(t\insl{r}t')\in\YoungT(\Si,\nu)$, with~$\r(t)=x_1\ldots x_k$ and~$\c(t)=y_1\ldots y_k$.
When computing~$(t\insl{r}t')$ (resp.~$(t\insl{c}t')$), one constructs a standard skew tableau $\Q(t\insl{r}~t')$ (resp.~$\Q(t\insl{c}t')$) of shape~$\nu/\mu$ over~$[n]$.
Start with~$\T^{k+1}=t'$ and an empty tableau~$\Q_r^{k+1}$ (resp.~$\Q_c^{k+1}$) of shape~$\mu$. For each~$i=k,\ldots,1$, compute~$x_i\insl{}\T^{i+1}$ (resp.~$y_i\insl{}\T^{i+1}$); let~$\T^i$ be the resulting super tableau. Add a
box filled with~$(k-i+1)$ to~$\Q_r^{i+1}$ (resp.~$\Q_c^{i+1}$) in the same place as the box that belongs to~$\T^i$ but not
to~$\T^{i+1}$; let~$\Q_r^i$ (resp.~$\Q_c^i$) be the resulting tableau. Output~$\Q_r^1$ (resp.~$\Q_c^1$) as $\Q(t\insl{r}t')$ (resp.~$\Q(t\insl{c}t')$).
Note that $\Q(t\insl{r}t')\neq \Q(t\insl{c}t')$, for some $t,t'$ in~$\YoungT(\Si)$.

When computing~$(t\insl{r}t')$ (resp.~$(t\insl{c}t')$), one constructs a skew tableau $\R(t\insl{r}t')$ (resp.~$\R(t\insl{c}t')$) of shape~$\nu/\mu$ over~$[n]$, as follows.
Start with~$\T^{k+1}=~t'$ and an empty tableau~$\R_r^{k+1}$ (resp.~$\R_c^{k+1}$) of shape~$\mu$. For each~$l=k,\ldots,1$, suppose that~$(i_l,j_l)$ is the position of~$x_l$ in~$t$ and compute~$x_l\insl{}\T^{l+1}$ (resp.~$y_l\insl{}\T^{l+1}$); let~$\T^l$ be the resulting super tableau. Add a
box filled with~$i_l$ to~$\R_r^{l+1}$ (resp.~$\R_c^{l+1}$) in the same place as the box that belongs to~$\T^l$ but not
to~$\T^{l+1}$; let~$\R_r^l$ (resp.~$\R_c^l$) be the resulting super skew tableau. Output~$\R_r^1$ (resp.~$\R_c^1$) as~$\R(t\insl{r}t')$ (resp.~$\R(t\insl{c}t')$).

\begin{lemma}
\label{L:Recordingtableau}
For all~$t$ and~$t'$ in~$\YoungT(\Si)$, we have~$\R(t\insl{r}t')=\R(t\insl{c}t')$.
\end{lemma}

Starting from~$\Q(t\insl{r}t')$ (resp.~$\Q(t\insl{c}t')$), we  recover~$\R(t\insl{r}t')$, as follows.
Start with an empty tableau~$\R^0$ of shape~$\mu$ and add a box filled with~$1$  in the same position as the box filled with~$1$ in~$\Q(t\insl{r}t')$ (resp.~$\Q(t\insl{c}t')$); let~$\R^1$ be the resulting skew tableau of filling~$\Tr^1$.
For each~$l=2,\ldots,k$, suppose that~$(i_l,j_l)$ is the position of~$l$ in~$\Q(t\insl{r}t')$ (resp.~$\Q(t\insl{c}t')$).
Add a box filled with~$\Tr^{l-1}(i_{l-1}, j_{l-1})$ if~$i_{l-1}\geq i_{l}$ or with~$\Tr^{l-1}(i_{l-1}, j_{l-1})+1$ otherwise (resp. with~$1$ if~$j_l>j_{l-1}$ or with~$\Tr(i_{l-1}, j_{l-1})+1$ otherwise) to~$\R^{l-1}$ in the position~$(i_l,j_l)$; let~$\R^l$ be the resulting skew tableau of filling~$\Tr^{l}$.
Output~$\R^k$ as $\R(t\insl{r}t')$.

Starting from~$\R(t\insl{r}t')$,  the skew tableau~$\Q(t\insl{r}t')$ is recovered as follows.
Start with an empty tableau~$\Q^0$ of shape~$\mu$.
Add~$\lambda_1$ boxes to~$\Q^0$ in the same positions as the boxes filled with~$1$ in~$\R(t\insl{r}t')$, and fill theses boxes by~$1$ to~$\lambda_i$, from bottom to top and from left to right; let~$\Q^1$ be the resulting skew tableau.
For each~$i=2,\ldots, h$, add~$\lambda_i$ boxes to~$\Q^{i-1}$ in the same positions as the boxes filled with~$i$ in~$\R(t\insl{r}t')$, and fill theses boxes by~$\lambda_{i-1}+1$ to~$\lambda_{1}+\ldots+\lambda_i$, from bottom to top and from left to right; let~$\Q^i$ be the resulting skew tableau.
Output~$\Q^l$ as~$\Q(t\insl{r}t')$.

Starting from~$\R(t\insl{r}t')$, one recovers~$\Q(t\insl{c}t')$, as follows.
Start with an empty tableau~$\Q^0$ of shape~$\mu$.
For each~$i=1,\ldots, h$, order the boxes filled with~$i$ in~$\R(t\insl{r}t')$ from the rightmost to the leftmost one, and suppose that the ordering is~$(i_1,j_1)_i,\ldots, (i_{\lambda_i},j_{\lambda_i})_i$.
Then order all the non-empty boxes in~$\R(t\insl{r}t')$ by setting that~$(i_m,j_m)_i$ precedes~$(i_{m'}, j_{m'})_{i'}$ if~$m'>m$ or~$m=m'$ and~$i'<i$. Relabel all the boxes and suppose that the ordering is~$(i_1,j_1),\ldots, (i_k,j_k)$. Add~$k$ boxes to~$\Q^0$ in the same positions as the non-empty boxes in~$\R(t\insl{r}t')$ and fill each box at the position~$(i_m, j_m)$ by~$(k-m+1)$ for~$m=1,\ldots, k$, with respect to the previous ordering.
Output the resulting standard skew tableau as~$\Q(t\insl{c}t')$.

As a consequence, we obtain the following result:

\begin{theorem}
\label{T:RecordingTableau}
There is a one-to-one correspondence between the set of skew tableaux $\R(t\insl{r}t')$ and the set of standard skew tableaux~$\Q(t\insl{r}t')$ (resp.~$\Q(t\insl{c}t')$), for all~$t,t'$ in~$\YoungT(\Si)$.
\end{theorem}

Denote by~$R_{\Si}$ the $\mathbb{Z}$-algebra  whose linear generators are the monomials in~$\P(\Si)$.
A generic element in~$R_{\Si}$ is realized by a formal sum of super plactic classes with coefficients from~$\mathbb{Z}$, and also by a formal sum of super tableaux. 
Construct a canonical homomorphism~$R_{\Si}\to\mathbb{Z}[X]$  sending each super  tableau~$t$ to its monomial~$x^t$ formed  by the product of~$x_i$, each occurring as many times in~$x^t$ as $i$ occurs in~$t$. 
Define~$S_\lambda$ (resp.~$S_{\lambda/\mu}$) in~$R_{\Si}$ to be the sum of all super  tableaux of shape~$\lambda$ (resp.~$\lambda/\mu$) over~$\Si$.
Consider~$\lambda, \mu, \nu$ in~$\Pr$ such that~$\mu, \lambda\subseteq \nu$. For any~$t''\in\YoungT(\Si,\nu)$, the \emph{super Littlewood--Richardson coefficient}, denoted by~$c_{\lambda,\mu}^{\nu}$, is the number of pairs~$(t, t')\in\YoungT(\Si,\lambda)\times\YoungT(\Si,\mu)$ such that~$(t\insl{r}t') =t''$. The following identities hold in~$R_{\Si}$,~\cite{Hage2021Super, Hage2022RSK}:
\[
S_\lambda S_\mu  =\underset{\nu}{\sum} c_{\lambda,\mu}^{\nu} S_\nu
\quad\text{ and }\quad
S_{\nu/\lambda}= \underset{\mu}{\sum} c_{\lambda,\mu}^{\nu}S_{\mu}.
\]

A word~$x_1\ldots x_k$ over~$[n]$ is \emph{Yamanouchi} if every left subword~$x_k \ldots x_l$ contains at least as many $1'$s as it does $2'$s, at least as many $2'$s as $3'$s, and so on for all positive integers.
A skew tableau over~$[n]$ has \emph{weight}~$ (\lambda_1, \ldots, \lambda_l)$ if its entries consist of $\lambda_1$ $1'$s, and so on up to $\lambda_l$ $l'$s.
A skew tableau~$S$ over~$[n]$ is a
\emph{Littlewood--Richardson tableau} if~$\r(S)$ is Yamanouchi.

Consider~$(t,t')\in \YoungT(\Si,\lambda)\times\YoungT(\Si,\mu)$.
By construction,~$\R(t\insl{r}t')$ is of weight~$\lambda$.
When computing~$\R(t\insl{c}t')$, each time we add a box filled with~$k$, we  previously added boxes filled with~$k-1,\ldots,  1$ in strictly higher rows, showing that~$\r(\R(t\insl{c}t'))$ is Yamanouchi.
We deduce by Lemma~\ref{L:Recordingtableau} that~$\R(t\insl{r}t')$ is a Littlewood--Richardson tableau of weight~$\lambda$.
Suppose now that we are given~$t''$ in~$\YoungT(\Si,\nu)$ such that~$(t\insl{r}t') =t''$, and we are given~$\R(t\insl{r}t')$. We compute~$\Q(t\insl{r}t')$, and then find~$t$. That which remains of~$t''$ will be~$t'$.
Given now any Littlewood--Richardson  tableau~$R$ of shape~$\nu/\mu$ and of weight~$\lambda$, for some partition~$\nu$ containing~$\mu$, and any~$t''$ in~$\YoungT(\Si,\nu)$, then there exist~$t$ in~$\YoungT(\Si,\lambda)$ and~$t'$ in~$\YoungT(\Si,\mu)$ such that~$(t\insl{r}t') =t''$ and~$\R(t\insl{r}t') =R$. Indeed, we compute~$\Q(t\insl{r}t')$, and find a word that corresponds to the row reading of a super tableau of shape~$\lambda$. We can also compute~$\Q(t\insl{c}t')$, and find a word that corresponds to the column reading of a super tableau of shape~$\lambda$. Note that the two constructions yield to the same super tableau~$t$ of shape~$\lambda$. That which remains of~$t''$ will be the super tableau~$t'$ of shape~$\mu$.
As a consequence, we deduce our main result:

\begin{theorem}
\label{T:LittlewoodRichardsonRule}
Let~$\lambda, \mu$ and~$\nu$ be in~$\Pr$ such that~$\mu\subseteq \nu$. The  coefficient~$c_{\lambda,\mu}^\nu$ is equal to the number of  Littlewood--Richardson  tableaux of shape~$\nu/\mu$  and of weight~$\lambda$.
\end{theorem}

Following~\cite{BonettiSenatoVenezia98}, we show all the results for any ordering of~$\Si$.
We recover the results of~\cite{Thomas78, LascouxSchutsenberger81} when~$\Si=\Si_1=~[n]$.
We recover the result obtained in~\cite{KangKwon01} using Kashiwara's theory of crystal bases for the  representations of the general linear Lie superalgebra, when~$\Si=\{1<\ldots<m<\overline{1}<\ldots <\overline{n}\}$.

\begin{small}
\renewcommand{\refname}{\Large\textsc{References}}
\bibliographystyle{plain}
\bibliography{biblioCURRENT}
\end{small}

\section*{Appendix 1}
\addcontentsline{toc}{section}{Appendix1}

\subsubsection*{Proof of Lemma~\ref{L:Recordingtableau}}
Proceed by induction on the number of rows of~$t$.
The result is obvious when~$t$ contains only one row, so suppose it is true for all super tableaux containing less than~$k$ rows.
Suppose now that~$t$ contains~$k$ rows. Denote by~$r_1$ its topmost row and by~$t_1$ the super tableau formed by its remaining rows. Then the following equality
\[
(t\insl{r}t') = (t_1\insl{r}(r_1\insl{r}t'))= (t_1\insl{c}(r_1\insl{r}t'))
\]
 holds in~$\YoungT(\Si)$. Since~$t_1$ contains~$k-1$ rows, we deduce by the induction hypothesis, that
\[
\R(t\insl{r}t')=\R(t_1\insl{c}(r_1\insl{r}t')).
\]
Now, prove that the recording tableau corresponding to
the insertion of~$r_1$ followed by the rightmost column of~$t_1$ followed by the rest of~$t_1$ by columns, is equal to the recording tableau corresponding to the insertion of the rightmost column of~$t$ followed by the topmost row of the remaining columns of~$t$ followed by the remaining columns of~$t_1$. Indeed, suppose that~$\r(r_1)=x_1\ldots x_k$ and the column reading of the rightmost column of~$t_1$ is~$y_1\ldots y_l$ with~$x_k=y_l$. Then, the equality
\[
\begin{array}{rl}
&(y_1\insl{r}(y_2\insl{r}\ldots(y_{l-1}\insl{r}(x_1\insl{r}(x_2\insl{r}\ldots(x_k\insl{r}t')\ldots)))\ldots))\\
&=
(x_1\insl{r}(x_2\insl{r}\ldots(x_{k-1}\insl{r}(y_1\insl{r}(y_2\insl{r}\ldots(y_l\insl{r}t')\ldots)))\ldots))
\end{array}
\]
holds in~$\YoungT(\Si)$. 
Indeed, one shows that the words~$y_1\ldots y_{l-1}x_1\ldots x_k$ and $x_1\ldots x_{k-1}y_1\ldots y_l$ are super plactic congruent, and then the equality is deduced by the cross-section property.
Following the column-bumping lemma~\cite[Proposition~2.10]{LaScalaNardozzaSenato06}, 
the box added by the insertion of~$x_i$ is in the same row or a higher row than the box added by the insertion of~$x_{i+1}$, for~$i=1,\ldots, k-1$, whereas
the box added by the insertion of~$y_i$ is in a lower row than the box
added by the insertion of~$y_{i+1}$,  for~$i=1,\ldots, k-1$,
showing the claim.
Hence, after inserting the rightmost column of~$t$, we interchange the order
of insertion of the topmost row of the remaining columns and the last but one
column of~$t$ in the same way as above without changing the recording tableau.
Successively, we interchange the order of insertion of the topmost row and the end
columns, we then insert~$t$ entirely by columns and the corresponding recording tableau remains unchanged.

\vfill
\begin{flushright}
\begin{small}
\noindent \textsc{Nohra Hage} \\
\url{nohra.hage@univ-catholille.fr} \\
ICL, Junia, Université Catholique de Lille, LITL, F-59000 Lille, France.
\end{small}
\end{flushright}

\vspace{0.25cm}

\begin{small}---\;\;\today\;\;-\;\;\hhmm\;\;---\end{small} \hfill
\end{document}

%% file: macros.tex
\newcommand{\periodafter}[1]{\ifstrempty{#1}{}{#1.}}
\titleformat{\section}[block]{\scshape\filcenter\LARGE}{\thesection.}{.5em}{}
\titleformat{\subsection}[block]{\bfseries\filcenter\large}{\thesubsection.}{.5em}{\medskip}
\titleformat{\subsubsection}[runin]{\bfseries}{\thesubsubsection.}{.5em}{\periodafter}%{}[.]
\titlespacing{\subsubsection}{0pt}{\topsep}{.5em}

\titleformat{\section}[block]{\scshape\filcenter\Large}{\thesection.}{.5em}{}
\titleformat{\subsection}[runin]{\bfseries}{\thesubsection.}{.5em}{}[.]
\titleformat{\subsubsection}[runin]{\bfseries}{\thesubsubsection.}{.5em}{}[.]
\titlespacing{\subsubsection}{0pt}{10pt}{.5em}

\theoremstyle{ntheorem}
  	\newtheorem{theorem}[subsection]{Theorem}
  	
	\newtheorem{lemma}[subsection]{Lemma}

\theoremstyle{definition}

	\makeatletter
\def\@equationname{equation}

\makeatother

\fancypagestyle{plain}{
\fancyhf{}\fancyfoot[LC,RC]{}
\fancyfoot[LE,RO]{$\thepage$}

}

\newcount\hh
\newcount\mm
\mm=\time
\hh=\time
\divide\hh by 60
\divide\mm by 60
\multiply\mm by 60
\mm=-\mm
\advance\mm by \time
\def\hhmm{\number\hh:\ifnum\mm<10{}0\fi\number\mm}

\UseTips
\SelectTips{eu}{11}

\newdir{ >}{{}*!/-10pt/@{>}}
\newdir{ -}{{}*!/-10pt/@{}}
\newdir{> }{{}*!/+10pt/@{>}}

\makeatletter

\xyletcsnamecsname@{dir4{}}{dir{}}
\xydefcsname@{dir4{-}}{\line@ \quadruple@\xydashh@}
\xydefcsname@{dir4{.}}{\point@ \quadruple@\xydashh@}
\xydefcsname@{dir4{~}}{\squiggle@ \quadruple@\xybsqlh@}
\xydefcsname@{dir4{>}}{\Tttip@}
\xydefcsname@{dir4{<}}{\reverseDirection@\Tttip@}

\xydef@\quadruple@#1{%
	\edef\Drop@@{%
		\dimen@=#1\relax
		\dimen@=.5\dimen@
		\A@=-\sinDirection\dimen@
		\B@=\cosDirection\dimen@
		\setboxz@h{%
			\setbox2=\hbox{\kern3\A@\raise3\B@\copy\z@}%
			\dp2=\z@ \ht2=\z@ \wd2=\z@ \box2
			\setbox2=\hbox{\kern\A@\raise\B@\copy\z@}%
			\dp2=\z@ \ht2=\z@ \wd2=\z@ \box2
			\setbox2=\hbox{\kern-\A@\raise-\B@\copy\z@}%
			\dp2=\z@ \ht2=\z@ \wd2=\z@ \box2
			\setbox2=\hbox{\kern-3\A@\raise-3\B@ \noexpand\boxz@}%
			\dp2=\z@ \ht2=\z@ \wd2=\z@ \box2
		}%
		\ht\z@=\z@ \dp\z@=\z@ \wd\z@=\z@ \noexpand\styledboxz@
	}%
}

\xydef@\Tttip@{\kern2pt \vrule height2pt depth2pt width\z@
	\Tttip@@ \kern2pt \egroup
	\U@c=0pt \D@c=0pt \L@c=0pt \R@c=0pt \Edge@c={\circleEdge}%
	\def\Leftness@{.5}\def\Upness@{.5}%
	\def\Drop@@{\styledboxz@}\def\Connect@@{\straight@{\dottedSpread@\jot}}}
	
\xydef@\Tttip@@{%
	\dimen@=.25\dimen@
 	\B@=\cosDirection\dimen@
	\setboxz@h\bgroup\reverseDirection@\line@ \wdz@=\z@ \ht\z@=\z@ \dp\z@=\z@
	{\vDirection@(1,-1)\xydashl@ \xyatipfont\char\DirectionChar}%
	{\vDirection@(1,+1)\xydashl@ \xybtipfont\char\DirectionChar}%
}

\xydef@\ar@form{
	\ifx \space@\next \expandafter\DN@\space{\xyFN@\ar@form}%
	\else\ifx ^\next \DN@ ^{\xyFN@\ar@style}\edef\arvariant@@{\string^}%
	\else\ifx _\next \DN@ _{\xyFN@\ar@style}\edef\arvariant@@{\string_}%
	\else\ifx 0\next \DN@ 0{\xyFN@\ar@style}\def\arvariant@@{0}%
	\else\ifx 1\next \DN@ 1{\xyFN@\ar@style}\def\arvariant@@{1}%
	\else\ifx 2\next \DN@ 2{\xyFN@\ar@style}\def\arvariant@@{2}%
	\else\ifx 3\next \DN@ 3{\xyFN@\ar@style}\def\arvariant@@{3}%
	\else\ifx 4\next \DN@ 4{\xyFN@\ar@style}\def\arvariant@@{4}%
	\else\ifx \bgroup\next \let\next@=\ar@style
	\else\ifx [\next \DN@[##1]{\ar@modifiers{[##1]}}%]
	\else\ifx *\next \DN@ *{\ar@modifiers}%
	\else\addLT@\ifx\next \let\next@=\ar@slide
	\else\ifx /\next \let\next@=\ar@curveslash
	\else\ifx (\next \let\next@=\ar@curveinout %)
	\else\addRQ@\ifx\next \addRQ@\DN@{\ar@curve@}%
	\else\addLQ@\ifx\next \addLQ@\DN@{\xyFN@\ar@curve}%
	\else\addDASH@\ifx\next \addDASH@\DN@{\defarstem@-\xyFN@\ar@}%
	\else\addEQ@\ifx\next \addEQ@\DN@{\def\arvariant@@{2}\defarstem@-\xyFN@\ar@}%
	\else\addDOT@\ifx\next \addDOT@\DN@{\defarstem@.\xyFN@\ar@}%
	\else\ifx :\next \DN@:{\def\arvariant@@{2}\defarstem@.\xyFN@\ar@}%
	\else\ifx ~\next \DN@~{\defarstem@~\xyFN@\ar@}%
	\else\ifx !\next \DN@!{\dasharstem@\xyFN@\ar@}%
	\else\ifx ?\next \DN@?{\ar@upsidedown\xyFN@\ar@}%
	\else \let\next@=\ar@error
	\fi\fi\fi\fi\fi\fi\fi\fi\fi\fi\fi\fi\fi\fi\fi\fi\fi\fi\fi\fi\fi\fi\fi \next@}

\makeatother

\newcommand{\qfl}{\xymatrix@1@C=10pt{\ar@4 [r] &}}

\DeclareMathOperator{\U}{\mathcal{U}}

\DeclareMathOperator{\Si}{\mathcal{S}}

\DeclareMathOperator{\YoungST}{Ys}
\DeclareMathOperator{\YoungT}{Yt}
\DeclareMathOperator{\SkewT}{St}

\renewcommand{\phi}{\varphi}
\renewcommand{\epsilon}{\varepsilon}

\newcommand{\M}{\mathbf{M}}

\renewcommand{\leq}{\leqslant}
\renewcommand{\geq}{\geqslant}

\newcommand{\insl}[1]{\rightsquigarrow_{#1}}
\newcommand{\insr}[1]{\;\raisebox{0.1em}{{\rotatebox[origin=c]{180}{$\rightsquigarrow$}}}_{#1}\;}

\definecolor{cyan}{RGB}{175,238,238} 
\definecolor{Red}{rgb}{0.96,0.17,0.20}
\definecolor{RedD}{rgb}{0.57, 0.0, 0.04}
\definecolor{Green}{rgb}{0.0,1.0,0.0}
\definecolor{GreenL}{rgb}{0.0,1.0,0.0} 
\definecolor{GreenD}{rgb}{0.64,0.76,0.68} 
\definecolor{Yellow}{rgb}{1.0,1.0,0.19}
\definecolor{YellowD}{rgb}{0.93, 0.84, 0.25}
\definecolor{Blue}{rgb}{0.0,1.0,1.0}
\definecolor{BlueD}{rgb}{0.63,0.79,0.95}
\definecolor{vert}{rgb}{0,0.45,0}
\definecolor{bazaar}{rgb}{0.6, 0.47, 0.48}
\definecolor{bronze}{rgb}{0.8, 0.5, 0.2}
\definecolor{darkspringgreen}{rgb}{0.09, 0.45, 0.27}

\def\Pr{\mathcal{P}}

\def\Tr{\mathcal{T}}
\def\Yr{\mathcal{Y}}
\makeatletter
\def\blfootnote{\xdef\@thefnmark{}\@footnotetext}
\makeatother
\def\r{\mathbf{r}}
\def\c{\mathbf{c}}
\def\R{\mathbf{R}}
\def\P{\mathbf{P}}
\def\T{\mathbf{T}}
\def\Q{\mathbf{Q}}

\newcommand{\ifthen}[2]{\ifthenelse{#1}{#2}{}}